\documentclass[12pt,thmsa]{article}
\usepackage{amsfonts}
\usepackage{amssymb}
\newtheorem{teo}{Theorem}[section]

\newtheorem{obs2}[teo]{Remark}

\newtheorem{tea}{Theorem}[subsection]

\newtheorem{no2}[teo]{Note}

\newtheorem{no3}[tea]{Note}

\newcommand{\Gal}{{\rm Gal}}
\newcommand{\Frob}{{\rm Frob }}

\newcommand{\mod}{{\rm mod}}

\newcommand{\Q}{\mathbb{Q}}

\newcommand{\PGL}{{\rm PGL}}
\newcommand{\PGSp}{{\rm PGSp}}

\newcommand{\GL}{{\rm GL}}

\newcommand{\F}{{\mathbb{F}}}

\newcommand{\PSp}{{\rm PSp}}

\newcommand{\GSp}{{\rm GSp}}
\newcommand{\Sp}{{\rm Sp}}
\begin{document}
\title{{\bf A non-solvable extension of $\Q$ unramified outside $7$
}
}
\author{Luis Dieulefait
\\
Dept. d'\'{A}lgebra i Geometria, Universitat de Barcelona;\\
Gran Via de les Corts Catalanes 585;
08007 - Barcelona; Spain.\\
e-mail: ldieulefait@ub.edu
%\thanks{Research partially supported     
%by MCYT grant .....}\\ 
  }
%MSC (2000): Primary 11F80, 12F12 
\date{\empty}

\maketitle

\vskip -20mm
%\titlerunning{Images of $3$-dimensional Galois representations
%

\begin{abstract}
We consider a mod $7$ Galois representation attached to a genus $2$ Siegel cuspforms of level $1$ and weight $28$ and using some of its Fourier coefficients and eigenvalues computed by N. Skoruppa and the classification of maximal subgroups of $\PGSp(4,p)$ we show that its image is as large as possible. This gives a realization of $\PGSp(4,7)$ as a Galois group over $\Q$ and the corresponding number field provides a non-solvable extension of $\Q$ which ramifies only at $7$.\\
\end{abstract}
%

%(Running head: .........)

\section{The result and its proof}

This brief note is aimed just to complement the remarkable work that Lassina Demb\'{e}l\'{e} and collaborators have done to solve Gross' conjecture about the existence of a non-solvable Galois extension of $\Q$ unramified outside $p$ for any prime $p$. The problem has an easy solution for primes $p>7$ by just inspecting the images of mod $p$ Galois representations attached to a few level $1$ classical cuspidal modular forms, and because of Serre's conjecture (in its strong form), it is easy to see that this approach can not solve the cases of small primes. Demb\'{e}l\'{e} and collaborators have solved the problem for $p=2,3,5$ (see [De] and [DGV]) by computing explicitly Hilbert modular forms over some totally real number field unramified outside $p$ with level also unramified outside $p$ and such that some of the corresponding residual representations in characteristic $p$ has non-solvable image. Their approach seems appropriate to attack also the case of $p=7$, it is just because of the cost of the required computations that they were not able to carry them on and exhibit an example solving Gross' conjecture for $p=7$.\\
The solution that we propose consists in looking at a mod $7$ representation attached to a level $1$ genus $2$ Siegel cuspform. The study of the images of these Galois representations was already developed in the author's thesis (see [Di1], chapter 6, and its published version [Di2], which contains less computations) where it was shown that the residual images were ``as large as possible" for almost every prime, provided that the form was not a Maass spezialform, and under certain irreducibility condition on one characteristic polynomial. In particular, for a Galois conjugacy class of cuspforms of level $1$ and weight $28$ it was shown that the image was generically large (there was also a previous result showing that, for this cuspform, the image was large for certain infinite set of inert primes containing $p=53$, see [DKR]).\\
Let us recall the setup from [Di1] (or [Di2]): we start from a genus $2$ Siegel cuspform $f$ of level $1$ and weight $28$, which is not a Maass spezialform (it is known to have multiplicity one). There is a unique conjugacy class, and the field $E$ generated by the eigenvalues of $f$ is the cubic field generated by some root of $P(x) = x^3 - x^2 - 294086 x - 59412960$. \\
Taylor and Weissauer proved the existence of a compatible family of $\lambda$-adic symplectic $4$-dimensional Galois representations $\{  \rho_\lambda \}$ attached to $f$ (see [T], [W]). This compatible system has conductor $1$, and the characteristic polynomial of the image of $\Frob \; p$ when $\lambda \nmid p$ is:
$$Pol_p(x) =  x^4 - a_p x^3 + (a_p^2 - a_{p^2} - p^{2k-4}) x^2 - a_p p^{2k-3} x + p^{4k-6}  $$
where $a_i$ denotes the $i$-th Hecke eigenvalue of $f$, for any $i$, and $k$ denotes the weight of $f$. \\
A priori, these $\lambda$-adic Galois representations are not known to be defined over $E_\lambda$, even if all characteristic polynomials have coefficients in $E$. Nevertheless, it was shown in [DKR] that one can work as if they were, i.e, just formally reduce mod $\lambda$ by reducing the coefficients of the characteristic polynomials and the resulting residual representation is known to be a quotient of $\rho_\lambda$. \\
We will consider a prime dividing $7$ in $E$ and the corresponding residual representation of $f$. The first Fourier coefficients of $f$ are available at [S], and have been quoted and used at [DKR] and [Di1]. Moreover, using these coefficients one can compute the first Hecke eigenvalues: $a_2, a_4, a_3, a_9, a_5, a_{25}$ (cf. [DKR] and [Di1]). Because we are only interested in a mod $7$ representation, let us just give the values of these eigenvalues in the residual representation. Let $\alpha$ be a root of $P(x)$. All Fourier coefficients and eigenvalues are computed in terms of $\alpha$. If we reduce $P(x)$ mod $7$ we obtain: 
$$ P(x ) \equiv (x+3) (x+4) (x+6) \quad \pmod{7} $$
Thus $7$ is split in $E$ and we can choose one of the three  residual representations with values on $\F_7$, and we choose to replace $\alpha$ by $1$ in order to compute, in $\F_7$, the residual value of the eigenvalues and the residual characteristic polynomials. This way we obtain:
$$a_2 = 4, a_4 = 5, a_3 = 3, a_9 = 2, a_5 = 1, a_{25} = 2 $$
And the following are the characteristic polynomials of the image of $\Frob \; p$ for $p=2,3,5$, factorized over $\F_7$:
$$Pol_2(x) \; \mod{7} = x^4 +  3 x^3 +   2 x^2  + 5  x +  2  $$
$$Pol_3(x) \; \mod{7} = (x+3) (x+4)  (x^2 + 4 x +5 )$$
$$Pol_5(x) \; \mod{7} =  (x^2 +    x  +      3) ( x^2 +         5  x  +    3) $$
Thus, the mod $7$ representation we are considering has image in $\GSp(4,7)$ containing three matrices with the above characteristic polynomials, and it cuts an extension of $\Q$ ramifying only at $7$. Because the residue field is prime (thus it has odd degree) and the multiplier being $\chi^{2k-3}$ it is known (cf. [DKR]) that if we consider the projectivization of this residual representation, if we see that its image contains $\PSp(4,7)$, then this image must be equal to $\PGSp(4,7)$. To show that the image is indeed large, we consider the classification of maximal subgroups of $\PGSp(4,q)$ given by Mitchell in geometric language and by Kleidman in group theoretic language (cf. [M] and [K], see also [KL1] and [Ki]). For the case of a prime field with $p = 7$  this classification reads as follows: \\
A  subgroup of $\PGSp(4,7)$, either contains $\PSp(4,7)$ or is contained in one of the following:\\
(1) a maximal parabolic subgroup (class $\mathcal{C}_1$ in [KL1]),\\
(2) a stabilizer of a decomposition of the underlying vector space in two two-dimensional subspaces (class $\mathcal{C}_2$ in [KL1]),\\
(3) a stabilizer of the extension field $\F_{7^2}$ (class $\mathcal{C}_3$ in [KL1]),\\
(4) a group isomorphic to $\PGL(2,7)$ (of type $\mathcal{S}$ in the notation of [KL1]),\\
(5) a group isomorphic to $ 2^4 \cdot O_4^{-}(2) \cdot 2$  (class $\mathcal{C}_6$ in [KL1]).\\
(6) a group isomorphic to $A_7 \cdot 2$ (of type $\mathcal{S}$ in the notation of [KL1]).\\

The image will be contained in a group as in (1) if and only if the representation is reducible. In case (2) there is a quadratic number field $K$ such that the restriction to $K$ of the representation is reducible, and for primes $p$ that are inert in $K$ and different from $7$ the trace of the representation at $p$ must be $0$ (see [Di1]). Case (3) is similar to the previous two cases, but after extending scalars to $\F_{p^2}$: over this extension a subgroup contained in a maximal subgroup in case (3) is again either reducible or contains a reducible normal subgroup of index $2$.\\
In case (5) the subgroup of order $16$ is an elementary abelian $2$-group, and the corresponding maximal subgroup of $\Sp(4,7)$ is known to be the normalizer of an extra-special group of order $32$ (cf. [K], [KL2]).\\

If we call a group as in (2) and (3) imprimitive, our goal is to show that the image of the mod $7$ representation is: \\
(a) absolutely irreducible, \\
(b)not imprimitive, and\\
(c) its projectivization is not contained in cases (4), (5) and (6).\\

(a) Assuming reducibility (over $\bar{\F}_{7}$), we can semisimplify for simplicity, and there will be either some $1$-dimensional or  some $2$-dimensional irreducible component. In the first case, this component corresponds to a character, but since the representation is unramified outside $7$, this character must be some power $\chi^i$ of the mod $7$ cyclotomic character. In particular, all characteristic polynomials must have at least one root in $\F_7$, contradicting the fact that both $Pol_2$ and $Pol_5$ don't have any such root. In the second case, the representation is the sum of two two-dimensional irreducible components (we can assume they are irreducible because the case of characters has just been taken care of). In this case, reducibility can not be taking place over $\F_7$ as in (1) because the image will be contained in two copies of $\GL(2,7)$ which is too small to contain an element with characteristic polynomial as $Pol_2$: an easy computation shows that a matrix with this characteristic polynomial has projective order $25$. On the other hand, if the image is only reducible after extending scalars, then the two irreducible components must be one conjugated to the other by the generator of $\Gal(\F_{49}/ \F_7)$. Observe also that the two components must have determinant defined over $\F_7$, because this determinant is a character of the absolute Galois group of $\Q$ unramified outside $7$. This is incompatible with the factorization of $Pol_3$, because the two components would have characteristic polynomials $(x+3) (x+4)$ and $  (x^2 + 4 x +5 )$ and since they have coefficients in $\F_7$ (and must be one conjugated to the other) they should be equal, but they are not.\\

(b) If the image were an (irreducible) imprimitive group, it should contain a normal reducible subgroup $H$ of index $2$, and if we consider the quadratic field $K$ fixed by $H$ for all primes inert in this field the trace should be $0$ (cf. [Di1]). Since the representation is unramified outside $7$, $K$ must be $\Q(\sqrt{-7})$ and in particular we should have $a_3 = 0$, which is not the case.\\

(c) Using $Pol_2$ we easily see that the projective image is not contained in cases (4), (5) or (6), because the projective order of this matrix is $25$ (the order of  $ O_4^{-}(2)$ is $120$). \\
%As for case (5), in this case the intersection of the image with $\PSp(4,7)$ is contained in $ 2^4 \cdot O_4^{-}(2)$. Let us restrict to this %intersection. Here the subgroup of order $16$ is an elementary abelian $2$-group, and the corresponding maximal subgroup of $\Sp(4,7)$ is known to be %the normalizer of an extra-special group of order $32$ (cf. [K], [KL2]). It follows from the study of these extra-special groups made in [G], Theorem %1, that the extension $ 2^4 \cdot O_4^{-}(2)$ is split. On the other hand, we consider the characteristic polynomial......... 
%Since $2$ is a square mod $7$ and the multiplier of the symplectic representation is $\chi^5$ the image of $\Frob \; 2$ falls in $\PSp(4,7)$.....

We conclude that the projectivization of the mod $7$ representation has image equal to $\PGSp(4,7)$, a non-solvable group. Let us record the result just proved in the following:

\begin{teo} Let $f$ be a genus $2$ Siegel cuspform of level $1$ and weight $28$ which is not a Maass spezialform. Consider the compatible family of symplectic $4$-dimensional Galois representations attached to $f$. Then, one of the residual representations in characteristic $7$ obtained from this family has projective image $\PGSp(4,7)$. In particular, this non-solvable group corresponds to an extension of $\Q$ ramifying only at $7$.
\end{teo}

\section{Bibliography}

[De] Demb\'{e}l\'{e}, L.; {\it
A non-solvable Galois extension of $\Q$ ramified at $2$ only},
C. R. Math. Acad. Sci. Paris {\bf 347} (2009) 111-116.\\

[DGV] Demb\'{e}l\'{e}, L.;   Greenberg, M.;  Voight, J.; {\it  Nonsolvable number fields ramified only at $3$ and $5$}, preprint, available at:  http://arxiv.org/abs/0906.4374\\

[DKR] Dettweiler, M.; Kuhn, U.;  Reiter, S., {\it On Galois Realizations via Siegel Modular Forms}, Math. Res. Lett. {\bf 8} (2001) 577-588.\\

[Di1] Dieulefait, L., {\it Modular Galois Realization of Linear Groups}, Ph. D. Thesis,  Universitat de Barcelona, 2001, \\ available at: wstein.org/people/dieulefait/lvdtesis.ps
 \\

[Di2] Dieulefait, L., {\it On the images of the Galois representations attached to genus $2$ Siegel modular forms}, J.  Reine Angew. Math.   {\bf 553} (2002)  183-200.\\

[Ki] King, O.; {\it The subgroup structure of finite classical groups in terms of geometric configurations}, in {\it Survey in Combinatorics, 2005}, B. S. Webb (Ed), Cambridge University Press (2005) 29-56.\\

%[G] Griess, R. Jr., {\it Automorphisms of extra special groups and nonvanishing degree $2$ cohomology}, Pac. J. Math {\bf 48} (1973) 403-422\\
[K] Kleidman, P. B., {\it The subgroup structure of some finite simple groups}, Ph. D. Thesis, 
Cambridge, 1986. \\

[KL1] Kleidman, P., Liebeck, M., {\it The subgroup structure of the finite classical groups}, London Math. Soc. LNS 129,   Cambridge University Press, 1990. \\

[KL2] Kleidman, P., Liebeck, M., {\it A survey of the maximal subgroups of the finite simple groups}, Geometriae Dedicata {\bf 25} (1988) 375-389\\

[M] Mitchell, H., {\it The subgroups of the quaternary abelian
linear group}, Trans. Amer. Math. Soc. {\bf 15} (1914) 379-396.\\

[S] Skoruppa, N., {\it Siegel Modular Forms of Genus $2$: Eigenforms of Level $1$}, tables available at: http://wotan.algebra.math.uni-siegen.de/$\sim$modi/
\\

[T] Taylor, R., {\it On the $\ell$-adic cohomology of Siegel
threefolds}, Invent. Math. {\bf 114} (1993) 289-310.\\

[W] Weissauer, R., {\it Four dimensional Galois representations},  in {\it Formes automorphes (II), Le cas du groupe $\GSp(4)$}, Ast\'{e}risque {\bf 302} (2005) 67-150.

\end{document}